\documentclass[12pt,a4paper]{amsart}

\usepackage{amsmath,amsthm,amssymb}
\theoremstyle{plain}
\newtheorem{thm}{Theorem}

\newtheorem{lema}{Lemma}
\newtheorem{prop}{Proposition}
\newtheorem{cor}{Corolary}

\theoremstyle{definition}
\newtheorem{defn}{Definition}

\theoremstyle{remark}

\newcommand{\IR}{\text{${\mathbb{R}}$}}
\newcommand{\IC}{\text{${\mathbb{C}}$}}
\newcommand{\IZ}{\text{${\mathbb{Z}}$}}

\newcommand{\IN}{\text{${\mathbb{N}}$}}

\newcommand{\IT}{\text{${\mathbb{T}}$}}

\begin{document}

\title[Generalized Fourier series ]{Generalized Fourier series and shift-invariant subspaces}
\thanks{\it Math Subject Classifications:
  42A24 (Primary);  41A30, 42A65, 42A8 (Secondary)}

\

\

\author{K. S. Kazarian}
\thanks{The work was supported in part by the research grant MTM 2010-15790; Dept. of Mathematics, Mod. 17, Universidad Aut\'onoma de Madrid, 28049, Madrid, Spain e-mail: kazaros.kazarian@uam.es}

%\markboth{Research Summary}{  Kazarian}

\begin{abstract}{A principal shift invariant subspace of $L^{2}(\IR)$ is isometric to a weighted norm space $L^{2}(\IT, w)$. Using results obtained earlier by the author on the basis properties of subsystems of the trigonometric system in the weighted norm spaces $L^{2}(\IT, w)$ we obtain corresponding results concerning generators of principal shift invariant subspaces.  }

\end{abstract}

\maketitle

%\vfill

 \vskip22pt

 \section{Introduction}

Shift invariant (or shortly SI) spaces  have been studied intensively in recent years  because of their importance in several areas
such as the theory of wavelets, spline systems, Gabor systems or approximation
theory.  They
are also used as models for spaces of
signals and images in mathematical and engineering applications. One of the main tools for the study of SI spaces is the Fourier analysis. In this context the
concept of generalized Fourier series (or GF series) can play an important role for obtaining constructive results when the generators of the SI are given. We
consider mainly univariate case because in the multivariate case the results on GF series are far from being complete.  Main literature on the shift-invariant
spaces can be found in the following monographs and  articles
[\cite{Ch:1}, \cite{D:92},\cite{HW:96},\cite{BDR:1}, \cite{BDR:2},\cite{JP:1},\cite{R:1}].

A closed subspace $V$ of $L^{2}(\IR)$ is called shift invariant if for any $f\in V$ and all $k\in \IZ$ $f(\cdot +k) \in V.$
If $\Gamma$ is a subset of $L^{2}(\IR)$ then we denote by $S(\Gamma)$ the shift invariant space
generated by $\Gamma$,
\[
S(\Gamma) =  \overline{\mbox{span}} \{\varphi(\cdot +k) : k\in \IZ, \varphi\in \Gamma \}.
\]
If $\Gamma$ consist of just  one element $\varphi$ then we will write $S(\varphi)$ instead of $S(\Gamma)$.

  We  consider that
the Fourier transform, an isometry from $L^2(\mathbb{R})$ onto itself, is defined so that the image
of a function $f\in L^1(\mathbb{R}) \cap
 L^2(\mathbb{R})$ is given by
 \[
  \widehat{f}({y}) = \int_{\mathbb{R}} f ({
 x}) e^{-2\pi i{x}\cdot {y}} d{x}, \qquad y\in \IR.
 \]

If $w\geq 0$ is a measurable function on a measurable set $E \subseteq \IR$  then we say that $\phi \in L^{2}(E,w)$ if $\phi:E\rightarrow \IC$ is measurable on $E$ and the norm is defined by
\[
 \| \phi \|_{L^{2}(E,w)}:= \left( \int_{E} |\phi(t)|^{2} w(t) dt\right)^{\frac{1}{2}} < +\infty.
 \]
A system  $ \{y_{k} \}_{k\in \IZ}$ in a Banach space $\bf{B}$ is called an $M-$basis if $ \{y_{k} \}_{k\in \IZ}$ is complete minimal in $\bf{B}$ and its dual
system is total with respect to the space $\bf{B}.$

  Let ${\IT} = {\IR}/{\IZ}.$ For $f,g\in L^2(\mathbb{R})$ following \cite{BDR:1} we put
 \[
 [{f}, {g}](t) = \sum_{{k}\in \IZ } {f}(t + k) {\overline{g(t + k)}}.
 \]
 It is easy to observe that $[{f}, {g}] \in L(\mathbb{\IT})$ and $|[{f}, {g}]|^{2} \leq [{f}, {f}] [{g}, {g}]$.
 When it will be  convenient we also use the notation $\Phi_{g} = [\widehat{g}, \widehat{g}]$ which is a well defined $1-$periodic non-negative function.

 Let $\varphi \in L^2(\mathbb{R})$ then
   by Plancherel's
theorem  we have that for any finite number of  coefficients $\{\gamma_{
k} \}_{|k|\leq N}$
\begin{align}\label{pl:1}
\int_{\IR}\bigg|\sum_{|k|\leq N} \gamma_{ k} \varphi(x+k)\bigg|^{2}dx & =
\int_{\IR}\bigg|\sum_{|k|\leq N} \gamma_{ k}e^{2\pi ik\cdot t}\bigg|^{2}
|\widehat{\varphi}({t})|^{2} dt\\ & = \int_{\IT} \bigg|\sum_{|k|\leq N} \gamma_{ k}e^{2\pi ik\cdot t}\bigg|^{2} \Phi_{\varphi}(t) dt.
\end{align}
Hence, one easily establishes an isometry between the
subspace $S(\varphi)$ and the subspace
\begin{equation}\label{sfin:1}
\widehat{S(\varphi)}= \widehat{\varphi} L^{2}(\IT, \Phi_{\varphi})
\subset L^{2}(\IR).
\end{equation}

If the function $\Phi_{\varphi}$ is not equal to zero on a set of positive measure we have that
$$
\frac{1}{\sqrt{\Phi_{\varphi}}} \widehat{\varphi} \in \widehat{S(\varphi)}.
$$
Hence the system $\{ \varphi_{0}(t+k)(x)\}_{k\in \IZ}$ where
\begin{equation}\label{or:1}
\widehat{\varphi_{0}}= \frac{1}{\sqrt{\Phi_{\varphi}}} \widehat{\varphi},
\end{equation}
  is an orthonormal basis of the subspace $S(\varphi)$(see \cite{D:92}, p. 140), \cite{S:94}, p. 20). Observe that $\varphi_{0}\in S(\varphi)$ by (\ref{sfin:1}).

 \

 The following claim clarifies if $\Phi_{\varphi}$ can be equal to any non negative integrable function defined on $\IT.$
\begin{prop}\label{pr:4} Let $\psi\in L^{2}(\IR)$ be such that $\| \psi\|_{2} =1$ and the system $\{ \psi(x +k)\}_{k\in \IZ}$ is orthogonal.
Then for any $h\in L^{2}(\IT),$  $h(x) = \sum_{k=-\infty}^{+\infty} c_{k} e^{2\pi ikx}$ we will have that
\[
\Phi_{\varphi}(t) = |h(t)|^{2},\quad \mbox{where} \quad \varphi(x) = \sum_{k=-\infty}^{+\infty} c_{k} \psi(x +k).
\]
\end{prop}

The above proposition can be checked easily using the following well known lemma
(cf.~\cite{D:92}, p. 132; ~\cite{HW:96})

\begin{lema}\label{lem:ap1}
The system $\{ \psi(\cdot +j): { j}\in
\IZ\},$ where $\psi\in L^2(\IR),$  is an orthonormal system if
and only if
\[
\sum_{{ j}\in \IZ}|\widehat{\psi}({ t}+  { j})|^2
 = 1 \qquad \mbox{for a.e.}\quad  { t}\in \IR  .
\]
\end{lema}
\

 \subsection{Shift invariant subspaces}

\

  The function $ \varphi$ could be smooth, band limited, symmetric etc., but
{\it a priori} the function $ \varphi_{0}$ may not have the same properties. Thus the study of conditions for which  the system $\{ \varphi(t+k)(x)\}_{k\in
\IZ}$
is a basis in some sense of the subspace $S(\varphi)$ has its own interest.

By equalities (\ref{pl:1})--(\ref{sfin:1}) it is easy to establish that  the metric
 properties of the system $\{ \varphi(t+k)(x)\}_{k\in \IZ}$ in $S(\varphi)$
 are the same as of the trigonometric system $\{ e^{2\pi itx}\}_{k\in \IZ}$ in the space
 $L^{2}_{\Phi_{\varphi}}(\IT).$ Let us give some of the simple consequences of the properties of the trigonometric system in weighted spaces.

\

 If for some $C>0$
\[
1/C \leq \Phi_{\varphi}(t) \leq C, \qquad \mbox{for a.e.} \quad t\in {\IT}
\]
then the  system $\{ e^{2\pi itx}\}_{k\in \IZ}$ is a Riesz basis in the weighted-norm space
$L^{2}_{\Phi_{\varphi}}(\IT)$ and consequently the system $\{ \varphi(t+k)(x)\}_{k\in \IZ}$ is a Riesz basis in $S(\varphi).$ If
\begin{equation}\label{w:1}
\frac{1}{\Phi_{\varphi}} \in L(\IT)
\end{equation}
it follows that $\{ \varphi(t+k)(x)\}_{k\in \IZ}$ is a  basis in its natural enumeration in $S(\varphi)$ if and only if the weight function $\Phi_{\varphi}$
belongs to the class $A_{2}$ (\cite{HMW:1}, \cite{K:01}). Recall the definition of the Muckenhoupt classes $A_{p}, p>1$ (\cite{Mu:1}). A non-negative function
$w\in L({\IT})$ belongs to the class $A_{p}, p>1$ if and only if for some $C_{p} >0$
\[
\frac{1}{|I|} \int_{I} w(x) dx \bigg[ \int_{I} w(x)^{-\frac{1}{p-1}} dx\bigg]^{p-1} \leq C_{p}
\]
for any interval $I\subset \IR$. In the case (\ref{w:1}) recall that the trigonometric system is an $M-$basis in $L^{2}_{\Phi_{\varphi}}(\IT)$ (\cite{K:1}).
 Hence, we have
\begin{prop}\label{pr:1}
Let (\ref{w:1}) holds. Then the system $\{ \varphi(t+k)(x)\}_{k\in \IZ}$ is an $M-$basis in $S(\varphi).$
\end{prop}
The dual of the system $\{ \varphi(t+k)(x)\}_{k\in \IZ}$ in $S(\varphi)$ has the following form $\{ g(t+k)(x)\}_{k\in \IZ},$ where $g \in S(\varphi)$ is
defined by
\[
\widehat{g} = \frac{1}{\Phi_{\varphi}} \widehat{\varphi}.
\]
Indeed, for $j,k \in \IZ$
\[
\begin{array}{cl}
\int_{\IR} \varphi(x+j) \overline{g(x+k)} dx = \int_{\IR} \widehat{\varphi}(t) \overline{\widehat{g}(t)}e^{2\pi i(j-k)t} dt = \\ \int_{\IT}\frac{1}{\Phi_{\varphi}}(t) \Phi_{\varphi}(t)e^{2\pi i(j-k)t} dt  = \delta_{jk},
\end{array} %\right.
\]
where $\delta_{jk}$ is Kronecker's delta function.

Let us study the case
\begin{equation}\label{w:2}
\frac{1}{\Phi_{\varphi}} \notin L(\IT).
\end{equation}
Recall some definitions from \cite{K:1}. We are going to use results only in the space $L^{2},$ hence we will modify the definitions for this case.  In the following definition we do not mention the degree of the singularity
because it will always be $2.$
\begin{defn} \label{Def:8} A function $w$, defined on the real axis, is said to have
a singularity   of order $k (k \in {\IN})$
at the point $x_{0} \in \IR$, if there is an interval $I$ with center in $x_{0}$  such
that $$\frac{(x - x_{0} )^{2(k-1)}}{w(x)}  \notin L(I)
 \quad \mbox{but} \quad \frac{(x - x_{0})^{2k}}{w(x)}  \in L(I). $$
\end{defn}

We are going to study the case when the weight function $\Phi_{\varphi}$ has singularities only in  finitely many points and of finite order.
Let $X = \{x_1, x_2, \dots , x_s\} \subset \IT $ be the set of all
points where the weight function $\Phi_{\varphi}$ has singularities and let $\Lambda = \{ \lambda_{j} \}_{j=1}^{s} \subset\IN$ be the set of corresponding
orders of those singularities. Given the pair $(X,\Lambda)$ one can define generalized Fourier series which generally speaking may have the same role in the
space $L^{2}_{\Phi_{\varphi}}(\IT)$ as the Fourier series in $L^{2}(\IT).$
Then by Theorem 1 of \cite{K:1} we obtain
\begin{prop}\label{pr:2} Let $\varphi \in L^{2}(\IR)$ be such that $\Phi_{\varphi}$ has singularities only in  finitely many points
$X = \{x_1, x_2, \dots , x_s\} \subset \IT $ and order of the singularity at the point $x_{j}$ is $\lambda_{j}$.
Let $\Omega \subset \IZ$ be a set of $|\Lambda|$
consecutive integers and $\Omega^{{c}}= \IZ\setminus \Omega$ and $|\Lambda| \geq 1.$ Then
the system $\{{\varphi}( t + k ) \}_{k\in \Omega^{{c}}}$ is an $M-$basis
 in $S(\varphi)$.
 \end{prop}

In fact we can say more about the system $\{{\varphi}( t + k )
\}_{k\in \Omega^{{c}}}.$  Let
\begin{equation}\label{du:1}
\widehat{\eta}_{ n }(t) =   \widehat{\varphi}(t)[\Phi_{\varphi}(t) ]\sp
{-1} [e\sp {i2\pi nt} - T_{n} (t)] \quad n\in \Omega^{c},
\end{equation}
 where
$T_{ n} (t) = \sum_{n\in \Omega } c_{ k}^{(n)} e^{i2\pi
kt}$ interpolates the function $e^{i2\pi nt}$ at every point
$x_{j} (1\leq j\leq s )$ up to the order $\lambda_{ j}$. Evidently,
$\widehat{\eta}_{ n} \in \widehat{S(\varphi)},$ hence $ {\eta}_{n} \in S(\varphi)$ and we
have that for any $k,n\in \Omega^{c}$
$$
\int_{\IR} \varphi(x+k) \overline{\eta}_{n}(x) dx
$$
$$= \int_{\IR}
\widehat{\varphi}(t) e\sp {i2\pi
kt}\overline{\widehat{\varphi}}(t)[\Phi_{\varphi}(t) ]\sp {-1} [e
^{-i2\pi nt} - \overline{T}_{n} (t)] dt = \delta_{kn}.
$$
 Thus the equations (\ref{du:1})
define the system biorthogonal to $\{{\varphi}( t + k ) \}_{k\in
\Omega^{{c}}}.$ By Theorem 1 of \cite{K:1} it follows that the system
\[
g_{n}(t) = [\Phi_{\varphi}(t) ]\sp {-1} [e
^{i2\pi nt} - {T}_{n} (t)]
\]
is complete in $L^{2}_{\Phi_{\varphi}}(\IT)$. For any $f\in S(\varphi)$ we have that $\widehat{f}(t) = \widehat{\varphi}(t) h(t),$ where $h \in
L^{2}_{\Phi_{\varphi}}(\IT).$ Hence, if
\[
\int_{\IR} f(x) \overline{\eta}_{n}(x) dx = 0, \quad \mbox{for all}\quad n\in \Omega^{c}
\]
then it follows that
\[
0= \int_{\IR} \widehat{\varphi}(t) h(t)\overline{\widehat{\varphi}}(t) [\Phi_{\varphi}(t) ]\sp {-1} [e
^{-i2\pi nt} - \overline{T}_{n} (t)] dt
\]
\[ = \int_{\IT} h(t) [e
^{-i2\pi nt} - \overline{T}_{n} (t)] dt
\]
for all $n\in \Omega^{c}.$
Thus $h(t)= 0$ a.e. on $\IT$ and
$\{\eta_{n}(t) \}_{n\in \Omega^{{c}}}$ is complete in $S(\varphi).$ Defining
$\widehat{\eta}_{ n}$ for $n\in \Omega$ in the same way as
 in (\ref{du:1}) one gets $ {\eta}_{ n}\equiv 0.$ Thus if we put
\begin{equation*}\label{co:1}
c_{n}(f) = \int_{\IR} f(x) \overline{\eta}_{n}(x) dx \qquad \mbox{for
any} \quad n\in \IZ,
\end{equation*}
then the coefficients of the functions $ {\varphi}( x + k ) (k\in
\Omega)$ will be equal to zero.

The corresponding result from \cite{K:0} gives
\begin{prop}\label{pr:3} Let $\Omega \subset \IZ$ be a set of $|\Lambda|$
consecutive integers. Then
the system $\{{\varphi}( t + k ) \}_{k\in \Omega^{{c}}}$ is not a  basis
 in $S(\varphi)$ for any enumeration.
 \end{prop}

Having in mind Propositions \ref{pr:3} and \ref{pr:2} it is natural to  study if the system
$\{{\varphi}( t + k ) \}_{k\in \Omega^{{c}}}$   is a summation basis in its natural enumeration. Our first option is the  $(C,\alpha)$ summation method.

Recall the definition of the $(C,\alpha)$ means of a series $\sum_{k=0}^{\infty} c_{k}$ (see \cite{Zy:59}, p. 130):
\[
\sigma_{n}^{\alpha} = \frac{S_{n}^{\alpha}}{A_{n}^{\alpha}} = \frac{1}{A_{n}^{\alpha}} \sum_{k =0}^{n} A_{n-k}^{\alpha} c_{k} =
\frac{1}{A_{n}^{\alpha}} \sum_{k =0}^{n} A_{n-k}^{\alpha-1} S_{k},
\]
where for $\alpha>0$ the numbers $A_{n}^{\alpha} (n=0,1,\dots)$ are determined by the following formula
\[
\sum_{n=0}^{\infty} A_{n}^{\alpha} x^{n} = (1-x)^{-\alpha-1}, \quad A_{n}^{\alpha} = \bigg(\begin{array}{cc} n&+\alpha \\ &n \end{array} \bigg) =
\frac{n^{\alpha}}{\Gamma(\alpha+1)},
\]
and $S_{n} = \sum_{k=0}^{n} c_{k}.$
The following result is a direct consequence of Theorem 3 from \cite{K:3}.
\begin{thm}\label{thm:1} The system
$\{{\varphi}( t + k ) \}_{k\in \Omega^{{c}}}$  in the natural enumeration will be a $(C,\alpha)$ summation basis in $S(\varphi)$ if and only if
the function $\Phi_{\varphi}$ can be represented simultaneously in the following forms
\begin{equation}\label{ca:1}
\Phi_{\varphi}(x)= \omega^{2}(x) w(x), \qquad  \Phi_{\varphi}(x)= \omega_{1}^{2}(x) w_{1}(x)
\end{equation}
where
\begin{equation}\label{cao:1}
\omega(x) = \prod_{j=1}^{s}\left|\sin \pi(x-x_{j})\right|^{\lambda_{j}}, \qquad  \omega_{1}(x) = \prod_{j=1}^{s}\left|\sin
\pi(x-x_{j})\right|^{\beta_{j}}
\end{equation}
$0<\beta_{j}\leq \alpha (1\leq j\leq s),$ and $w, w_{1}$ are nonnegative functions from the class $A_{2}.$
\end{thm}

\

\subsection{Compactly supported generators}

Let $\varphi\in L^{2}(\IR)$ be a compactly supported complex valued
function.  We  suppose that the support of a function $f \in
L^{2}(\IR)$ is defined only modulo a null set as  $\{x\in \IR:
f(x)\neq 0\}.$ It is well known that the Fourier coefficients of the function
$\Phi_{\varphi}$ satisfy to the following equations:
\begin{equation*}\label{fc:1}
c_{k}(\Phi_{\varphi}) = \int_{0}^{1}\Phi_{\varphi}(t)e^{-2\pi
ikt}dt=\int_{\IR}\varphi(t)\varphi(t+k)dt \qquad \forall k\in \IZ.
\end{equation*}
Hence, $\Phi_{\varphi}$ is a trigonometric polynomial. Let $X = \{x_1, x_2, \dots , x_s\} \subset \IT $ be the set of all
points where the trigonometric polynomial $\Phi_{\varphi}$ has a
zero. We will
consider the case when $X \neq \emptyset,$ which can be studied
using the theory of generalized Fourier series. In this case we have
\begin{equation*}\label{eq:1}
\Phi_{\varphi}(t) = |\omega(t)|^2P(t),
\end{equation*}
where $P(t)>0$ is a trigonometric polynomial,
\begin{equation*}\label{eq:2}
\omega(x)= \omega_{X,\Lambda}(x)
:= \prod_{j=1}^s \big[\sin \frac{x-x_j}{2}
\big]^{\lambda_{j}}
\end{equation*}
and  $\Lambda= \{ \lambda_{j}
\}_{l=1}^{s} \subset \IN.$
If  $\Omega \subset \IZ$ is a set of $|\Lambda|$
consecutive integers
then  Propositions \ref{pr:2}, \ref{pr:3} and Theorem \ref{thm:1} can be applied.
Observe that the equations (\ref{du:1})
define the system biorthogonal to $\{{\varphi}( t + k ) \}_{k\in
\Omega^{{c}}}.$ From Theorem 2 of \cite{K:1} (see also \cite{K:2}) follows
\begin{thm}\label{thm:2} For any $f\in S(\varphi)$
\[
\lim \sb {r \to {1-}} ||f - \sum \sb {n=-\infty } \sp {+\infty } r\sp {|n|}
c\sb n (f) {\varphi}( \cdot + n ) ||_{L^{2}(\IR)} = 0,
\]
where
\[
c_{n}(f) = \int_{\IR} f(x) \overline{g}_{n}(x) dx \qquad \mbox{for
any} \quad n\in \IZ\setminus \Omega.
\]
\end{thm}
Set
$\lambda_{0}:= \max_{1\leq j\leq s} \lambda_{j}.$ By Theorem \ref{thm:1}
 we obtain

\begin{cor}\label{cor:3} For any $f\in S(\varphi)$
$$\lim \sb {n \to +\infty} ||f(\cdot) - S_{n}^\alpha(f,\cdot)||_{L^{2}(\IR)} = 0$$
where
$$S_{n}^\alpha(f,x)=(1/A_n^\alpha) \sum_{k=-n}^{n} A_{n-|\nu |}^{\alpha-1}c_{k}(f){\varphi}( x + k ),$$
if and only if $\alpha > \lambda_{0}-\frac{1}{2}.$
\end{cor}

\

\subsection{An example}

\

 Let $\phi(t) = \chi_{[0,1]}(t).$ Then we have that
\[
\widehat{\phi}(\xi) = \int_{\IR} \phi(t) e^{-2\pi i\xi t} dt = \int_{[0,1]} e^{-2\pi i\xi t} dt = -\frac{1}{2\pi i\xi} [e^{-2\pi i\xi } -1]
\]
\[
=\frac{\sin {\pi \xi} \cos{\pi \xi} }{\pi \xi} + \frac{\sin^{2} \pi \xi }{\pi i\xi} = \frac{\sin {\pi \xi} }{\pi \xi} e^{-\pi i\xi }.
\]
It is evident that the shifts of $\phi$ constitute an orthonormal system. Hence, by Lemma \ref{lem:ap1} we will have
\begin{equation*}
\Phi_{\phi}(\xi) = \frac{|\sin {\pi \xi} |^{2}}{\pi^{2}}\sum_{k\in \IZ} \frac{1}{| \xi +k|^{2}} = 1 \qquad \mbox{a.e. on}\, \IR.
\end{equation*}
Which yields the following identity
\begin{equation}\label{eq:kar1}
 \sum_{k\in \IZ} \frac{1}{| \xi +k|^{2}} = \frac{\pi^{2}}{|\sin {\pi \xi} |^{2}} \qquad \mbox{a.e. on}\, \IR.
\end{equation}

\

 Let $\phi_{(1)}(t) = \chi_{[0,1]}\ast\chi_{[0,1]}(t).$ Then
\[
\widehat{\phi_{(1)}}(\xi) = \widehat{\phi}(\xi) \widehat{\phi}(\xi) = \frac{\sin^{2} {\pi \xi} }{\pi^{2} \xi^{2}} e^{-2\pi i\xi }.
\]
Afterwards, we have to calculate
\[
\Phi_{\phi_{(1)}}(\xi) = \frac{|\sin {\pi \xi} |^{4}}{\pi^{4}}\sum_{k\in \IZ} \frac{1}{( \xi +k)^{4}} = \frac{|\sin {\pi \xi} |^{4}}{6\pi^{2}} \frac{d^{2}}{d\xi^{2}} \frac{1}{|\sin {\pi \xi} |^{2}},
\]
where the equalities should be understood a.e. on $\IR$. We skip the details of the proof of the last equation because it is elemental.

Thus we obtain that
\[
\Phi_{\phi_{(1)}}(\xi) =  \frac{1}{3} [2 + \cos 2\pi \xi].
\]
We have also obtained the following identity
\begin{equation}\label{eq:kar2}
 \sum_{k\in \IZ} \frac{1}{| \xi +k|^{4}} = \frac{{\pi^{4}}}{3|\sin {\pi \xi} |^{4}} [2 + \cos 2\pi \xi] \qquad \mbox{a.e. on}\, \IR.
\end{equation}

Let $\psi_{(1)}(t) = \phi_{(1)}(2t) - \phi_{(1)}(2t-2).$ Then
\[
\widehat{\psi_{(1)}}(\xi) = \frac{1}{2} \widehat{\phi_{(1)}}(\frac{1}{2} \xi) - \frac{1}{2} \widehat{\phi_{(1)}}(\frac{1}{2} \xi) e^{-2\pi i\xi } = \frac{1}{2} \widehat{\phi_{(1)}}(\frac{1}{2} \xi) [1- e^{-2\pi i\xi }]
\]
\[
= \widehat{\phi_{(1)}}(\frac{1}{2}\xi ) i \sin (\pi \xi) e^{-\pi i\xi } = i 4  \frac{\sin^{2} {\pi \frac{\xi}{2}} }{\pi^{2} \xi^{2}} e^{-\pi i\xi }
\sin (\pi \xi) e^{-\pi i\xi }
\]
\[
 = i 4  \frac{\sin^{2} {\pi \frac{\xi}{2}} }{\pi^{2} \xi^{2}}  \sin (\pi \xi) e^{-2\pi i\xi }
\]
Then we will have
\[
\Phi_{\psi_{(1)}}(\xi) = \frac{4|\sin {\pi \xi} |^{2}}{\pi^{4}}\sum_{k\in \IZ} \frac{[1- \cos({\pi (\xi +k)})]^{2}}{( \xi +k)^{4}}.
\]
Afterwards we have to use the identity (\ref{eq:kar2})
\[
\sum_{k\in \IZ} \frac{[1- \cos({\pi (\xi +k)})]^{2}}{( \xi +k)^{4}} = \sum_{k\in 2\IZ +1} + \sum_{k\in 2\IZ}
\]
\[
= \frac{1}{16}[1+ \cos({\pi \xi })]^{2}\sum_{m\in \IZ} \frac{1}{( \frac{\xi+1}{2} + m)^{4}} + \frac{1}{16}[1- \cos({\pi \xi })]^{2}\sum_{m\in \IZ} \frac{1}{( \frac{\xi}{2} + m)^{4}}
\]
\[
= \frac{\pi^{4}}{48}[1+ \cos({\pi \xi })]^{2} \frac{[2 - \cos (\pi \xi)]}{|\sin {\pi \frac{\xi+1}{2}} |^{4}}  + \frac{\pi^{4}}{48}[1- \cos({\pi \xi })]^{2} \frac{[2 + \cos (\pi \xi)]}{|\sin {\pi \frac{\xi}{2}} |^{4}}
\]
\[
= \frac{\pi^{4}}{3}[1+ \cos({\pi \xi })]^{2} \frac{[2 - \cos (\pi \xi)]}{[1+ \cos({\pi \xi }) ]^{2}}
 \]
 \[ + \frac{\pi^{4}}{3}[1- \cos({\pi \xi })]^{2} \frac{[2 + \cos (\pi \xi)]}{[1- \cos({\pi \xi }) ]^{2}} = \frac{4\pi^{4}}{3}.
\]
Thus we obtain that
\[
\Phi_{\psi_{(1)}}(\xi) = \frac{16 }{3} \sin^{2} {\pi \xi}
\]
Hence, in $\IT$ $\Phi_{\psi_{(1)}}$ has singularity only  at the point $0$, a singularity of order $1$ at the point $0$.
By Theorem \ref{thm:1} we obtain the following
\begin{cor} For any $\alpha > \frac{1}{2}$ the system $\{ \psi_{(1)}(x +k) \}_{k\in \IZ\setminus\{0\}}$ in its natural enumeration is a $(C,\alpha)$ summation basis.
\end{cor}

\

\end{document}